\documentclass[11pt,letter,leqno]{amsart}
\usepackage{amssymb}
\numberwithin{equation}{section}

\newtheoremstyle{theorem}{3pt}{3pt}%
{\it}
{}
{\bfseries}
{:}
{.5em}
{}
\theoremstyle{theorem}
\newtheorem{theorem}{Theorem}[section]

\newtheorem{proposition}[theorem]{Proposition}

\newtheorem{lemma}[theorem]{Lemma}
\newtheorem{definition}[theorem]{Definition}

\newtheorem{definition-proposition}{Proposition/Definition}

\newtheoremstyle{example}{3pt}{3pt}%
{}
{}
{\sc}
{:}
{.5em}
{}
\theoremstyle{example}
\newtheorem{example}[theorem]{Example}
\newtheoremstyle{remark}{3pt}{3pt}%
{}
{}
{\sc}
{:}
{.5em}
{}
\theoremstyle{remark}

\numberwithin{equation}{section}

\newcommand{\thismonth}{\ifcase\month\or
  January\or February\or March\or April\or May\or June\or
  July\or August\or September\or October\or November\or December\fi
  \space\number\year}
\makeatletter
\newcommand{\low}{\@ifnextchar^{}{^{\vphantom x}}}
\newcommand{\high}{\@ifnextchar_{}{_{\vphantom I}}}
\makeatother
\DeclareSymbolFont{script}{U}{eus}{m}{n}
\DeclareSymbolFontAlphabet{\mathscr}{script}
\DeclareMathSymbol{\EuWedge}{0}{script}{"5E}
\DeclareMathAlphabet{\mathrmsl}{OT1}{cmr}{m}{sl}
\newcommand{\rssymb}[2]{\newcommand{#1}{{\mathrmsl{#2}}}}
\newcommand{\calsymb}[2]{\newcommand{#1}{{\mathcal{#2}}}}
\newcommand{\bbsymb}[2]{\newcommand{#1}{{\mathbb{#2}}}}
\newcommand{\lieoper}[2]{\newcommand{#1}{\mathop
  {\mathfrak{#2}\null}\nolimits}}
\newcommand{\oper}[3][n]{\newcommand{#2}{\mathop
  {\mathrm{#3}\null}\ifx n#1\nolimits\else\limits\fi}}
\newcommand{\rsoper}[3][n]{\newcommand{#2}{\mathop
  {\mathrmsl{#3}\null}\ifx n#1\nolimits\else\limits\fi}}
\bbsymb\C{C} \bbsymb\F{F} \bbsymb\HQ{H}\bbsymb\N{N} \bbsymb\Q{Q}
\bbsymb\R{R} \bbsymb\U{U} \bbsymb\V{V} \bbsymb\W{W} \bbsymb\Z{Z}
\bbsymb\bbf{F} \bbsymb\bbk{K} \bbsymb\bbi{I} \bbsymb\bbl{L} \bbsymb\bbo{O}
\bbsymb\bbj{J}
\bbsymb\bby{Y}
\bbsymb\bbp{P}
\bbsymb\bba{A}
\calsymb\cA{A} \calsymb\cB{B} \calsymb\cC{C} \calsymb\cD{D} \calsymb\cE{E}
\calsymb\cF{F} \calsymb\cG{G} \calsymb\cH{H} \calsymb\cI{I} \calsymb\cJ{J}
\calsymb\cK{K} \calsymb\cL{L} \calsymb\cM{M} \calsymb\cN{N} \calsymb\cO{O}
\calsymb\cP{P} \calsymb\cQ{Q} \calsymb\cR{R} \calsymb\cS{S} \calsymb\cT{T}
\calsymb\cU{U} \calsymb\cV{V} \calsymb\cW{W} \calsymb\cX{X} \calsymb\cY{Y}
\calsymb\cZ{Z}
\renewcommand{\geq}{\geqslant} \renewcommand{\leq}{\leqslant}
\oper\End{End} \oper\Hom{Hom}                    
\oper\Sym{Sym} \oper\Skew{Skew}
\oper\Aut{Aut}                                   
\oper\GL{GL} \oper\SL{SL}\oper\Symp{Sp}
\oper\CO{CO} \oper\On{O} \oper\SO{SO} \oper\Pin{Pin} \oper\Spin{Spin}
\oper\CU{CU} \oper\Un{U} \oper\SU{SU}
\rsoper\Diff{Diff} \rsoper\SDiff{SDiff}
\lieoper\der{der}                                
\lieoper\gl{gl} \lieoper\sgl{sl}\lieoper\symp{sp}
\lieoper\co{co} \lieoper\so{so} \lieoper\spin{spin}
\lieoper\cu{cu} \lieoper\un{u}  \lieoper\su{su}
\rsoper\Vect{Vect} \rsoper\Ham{Ham}
\newcommand{\norm}[2][]{|\mkern-2mu|#2|\mkern-2mu|
  _{\lower1pt\hbox{${}_{#1}$}}}
\newcommand{\Norm}[2][]{\bigl|\mkern-3mu\bigr|#2\bigr|\mkern-3mu\bigr|
  _{\lower1pt\hbox{${}_{#1}$}}}

\rsoper\dimn{dim}                           
\rsoper\grad{grad}                          
\rsoper\kernel{ker}\rsoper\image{im}        
\rsoper\alt{alt}   \rsoper\sym{sym}         
\rsoper\Ad{Ad}     \rsoper\ad{ad}           
\rsoper\CoAd{CoAd} \rsoper\coad{coad}       
\rsoper\trace{tr}  \rsoper\trfree{tf}       
\rsoper\detm{det}                           
\rsoper\Vol{Vol}                            
\rsoper\divg{div}                           
\rsoper\sign{sign}                          
\rssymb\iden{id}                            
\rssymb\vol{vol}                            
\oper\Imag{Im}\oper\Real{Re}                
\newcommand{\sd}{{\raise1pt\hbox{$\scriptscriptstyle +$}}}
\newcommand{\asd}{{\raise1pt\hbox{$\scriptscriptstyle -$}}}
\newcommand{\sdasd}{{\raise1pt\hbox{$\scriptscriptstyle\pm$}}}
\newcommand{\asdsd}{{\raise1pt\hbox{$\scriptscriptstyle\mp$}}}

\rsoper\scal{scal}
\def\kahl/{k\"ahler}
\def\Kahl/{K{\"a}hler}
\newcommand{\bfw}{{\mathbf w}}
\newcommand{\bfx}{{\mathbf x}}

\newcommand{\bfa}{{\mathbf a}}

\def\decdnar#1{\phantom{\hbox{$\scriptstyle{#1}$}}
\left\downarrow\vbox{\vskip15pt\hbox{$\scriptstyle{#1}$}}\right.}
\def\bba{{\mathbb A}}

\def\bbc{{\mathbb C}}

\def\bbf{{\mathbb F}}

\def\bbi{{\mathbb I}}
\def\bbj{{\mathbb J}}
\def\bbk{{\mathbb K}}
\def\bbl{{\mathbb L}}

\def\bbo{{\mathbb O}}
\def\bbp{{\mathbb P}}
\def\bbq{{\mathbb Q}}
\def\bbr{{\mathbb R}}

\def\bby{{\mathbb Y}}
\def\bbz{{\mathbb Z}}

\def\gra{\alpha}

\def\grg{\gamma}

\def\grl{\lambda}

\def\gro{\omega}

\def\grt{\tau}

\def\grz{\zeta}
\def\grD{\Delta}

\def\grG{\Gamma}

\def\grL{\Lambda}

\def\grS{\Sigma}

\def\calu{{\mathcal U}}
\def\cald{{\mathcal D}}

\def\calf{{\mathcal F}}

\def\cals{{\mathcal S}}

\def\calu{{\mathcal U}}

\def\calz{{\mathcal Z}}

\def\Se{Sasakian-Einstein }
\def\Ke{K\"ahler-Einstein }
\def\la#1{\hbox to #1pc{\leftarrowfill}}
\def\ra#1{\hbox to #1pc{\rightarrowfill}}
\def\Ric{{\rm Ric}}
\def\tU{{\tilde U}}

\def\fract#1#2{\raise4pt\hbox{$ #1 \atop #2 $}}
\def\decdnar#1{\phantom{\hbox{$\scriptstyle{#1}$}}
\left\downarrow\vbox{\vskip15pt\hbox{$\scriptstyle{#1}$}}\right.}

\begin{document}

\title{Sasakian Geometry and Einstein Metrics on Spheres}

\author{Charles P. Boyer and Krzysztof Galicki}
\address{Department of Mathematics and Statistics,
University of New Mexico,
Albuquerque, NM 87131.}
\email{cboyer@math.unm.edu}

\email{galicki@math.unm.edu}


\maketitle

\section{Introduction}

This paper is based on a talk presented by the first author at the
Short Program on Riemannian Geometry that took place at the Centre
de Recherche Math\'ematiques, Universit\'e de Montr\'eal, during
the period June 28-July 16, 2004. It is a report on our joint work
with J\'anos Koll\'ar \cite{BGK03} concerning the existence of an
abundance of Einstein metrics on odd dimensional spheres,
including exotic spheres. The article appeared electronically in
September of 2003, and answered in the affirmative the heretofore
open question of the existence of Einstein metrics on exotic
spheres. Evan Thomas helped us with the computer programs
\cite{BGKT03}.

\subsection{A Brief History of our Approach}
\medskip

For more than ten years the authors and their collaborators have
established a program employing Sasakian geometry to prove the
existence of Einstein metrics on compact odd dimensional
manifolds. It began \cite{BGM93a, BGM93b, BGM94a} with the study
of 3-Sasakian manifolds which are automatically Einstein. The main
technique used to establish the existence of Einstein 3-Sasakian
metrics was that of symmetry reduction, a method well-known in
symplectic geometry. This work reached its pinnacle in the {\it
Inventiones} paper \cite{BGMR}, where  we constructed all toric
3-Sasakian manifolds in dimension 7. In particular, we gave the
first examples of Einstein metrics on compact simply-connected
manifolds with arbitrary second Betti number. A survey of the
results in this area appeared later in \cite{BG99}.

In \cite{BG00a} we began the study of \Se structures. In this more
general case the symmetry reduction yielded very few results, and
another method was needed. In late 1999 we noticed the preprint
version of the paper \cite{DeKo01} by Demailly and Koll\'ar on the
LANL ArXivs which gave sufficient algebraic conditions for the
vanishing of the obstructions to solving the Monge-Amp\`ere
equations for K\"ahler orbifolds with positive first Chern class.
More importantly, Demailly and Koll\'ar exhibited explicit
examples of log del Pezzo surfaces with orbifold singularities as
certain hypersurfaces in weighted projective spaces. This entailed
showing that the singularities of the pair (log del Pezzo surface,
anticanonical divisor) is relatively mild, a condition known in
the algebraic geometry literature as {\it Kawamata log terminal},
or {\it klt} for short. Their article prompted us to develop a
method \cite{BG01b} which uses the links of hypersurface
singularities of weighted homogeneous polynomials to give the
first non-regular examples of \Se metrics in dimension five.
Shortly thereafter we received a preliminary version of \cite{jk1}
in which many more explicit examples of log del Pezzo surfaces
with a K\"ahler-Einstein orbifold metric were given. We teamed up
with our algebraic geometer colleague Michael Nakamaye
\cite{BGN03c,BGN02b,BG03} to apply the method of \cite{BG01b} to
prove the existence of many \Se metrics on the connected sums
$\#k(S^2\times S^3)$ for $1\leq k\leq 9.$ (And very recently Koll\'ar \cite{Kol04} has proven 
that \Se metrics exist on $k\#(S^2\times S^3)$ for all
$k\geq 6.$ In fact, he has shown that there are infinitely many
$2(k-1)$-dimensional families of \Se metrics on $k\#(S^2\times
S^3).$ Thus, all connected sums of $S^2\times S^3$ admit many \Se
metrics.)  However, the question of
how to handle the klt conditions when the K\"ahler orbifolds have
branch divisors alluded us. This occurs, for example, when the
links are homotopy spheres, and we had early on recognized that it
was a potential approach to proving the existence of an Einstein
metric on an exotic sphere.  In fact, we already had what we
thought was a good candidate for a \Se structure on an exotic
sphere. Then in the summer of 2003 we began serious discussions
with J\'anos Koll\'ar who understood the singularities with branch
divisors present. We not only were able to verify that our
original candidate provided the first example of an Einstein
metric on an exotic sphere, but we found, much to our surprise and
delight, that Einstein metrics exist in huge numbers on many
homotopy spheres. We would like to take this opportunity to thank
J\'anos for taking an early interest in our work, and then
providing much of the essential ingredients of our joint papers. We also thank him for useful 
comments concerning this expository paper.

\subsection{A History of Einstein Metrics on Spheres}

The name comes from Einstein famous work creating general
relatively. Working with the indefinite Lorenzian signature of
4-dimensional space-time, Einstein's reasoning went roughly as
follows: the total amount of energy and momentum in the universe
should equal the curvature of the universe. Energy and momentum is
represented by a symmetric 2-tensor $T_{\mu\nu}$, and there are
precisely two natural symmetric 2-tensors in Riemannian
(Lorenzian) geometry, the Ricci curvature, $R_{\mu\nu}$, and the Riemannian
metric itself $g_{\mu\nu}.$ Hence, in tensor indices one has the equation
$$G_{\mu\nu}:=R_{\mu\nu}-\frac{1}{2}s g_{\mu\nu}=\kappa T_{\mu\nu}-\Lambda 
g_{\mu\nu},$$
where  $G_{\mu\nu}$ is sometimes called the Einstein tensor,
$s=\sum_{\mu,\nu}g^{\mu\nu}R_{\mu\nu}$ is the scalar curvature,
and $\kappa,\Lambda$ are constants. If one assumes that the only
energy in the universe is gravitational, then $T_{\mu\nu}=0.$
Actually, Einstein originally had  the ``cosmological constant''
$\Lambda=0$ as well, but he then inserted a non-zero value of
$\Lambda$ to obtain a solution representing a static universe.
Later when Hubble discovered that the universe was actually
expanding, he called inserting $\Lambda$ ``the biggest blunder of
my life''. Ironically, the recently discovered acceleration in the
expansion of the universe again suggests a non-zero value of
$\Lambda$ in order to account for the so-called ``dark energy''
which represents a type of anti-gravitational pressure that causes
the acceleration in the expansion of the universe. So in the absence of other forces, the 
Einstein equations become $R_{\mu\nu}=(\frac{1}{2}s-\grL) g_{\mu\nu}.$ This is the origin of 
the mathematicians' well known definition of an Einstein metric, viz.

\begin{definition}\label{Emetric} A Riemannian metric $g$ is said to 
be an {\bf Einstein metric} if
$\Ric_g=\grl g$ for some constant $\grl.$
\end{definition}

The fact that $\grl$ is a constant is a consequence of the Bianchi
identities. The best known examples of Einstein metrics are
metrics of constant curvature, and long before Einstein gave us
his modern view of the universe, even before Riemann's
``epoch-making'' essay on the foundation of geometry, Gauss had
studied spherical geometry and understood the ``round sphere
metric'' at least in dimension 2. More than one hundred years had
passed since the deaths of both Gauss and Riemann before an
example of an Einstein metric other than the round sphere metric
was shown to exist on a sphere. In 1973, Gary Jensen \cite{Jen73}
proved the existence of what later became known as a ``squashed''
metric on $S^{4n+3}.$ The next 25 years brought about very little
change with just a handful of new results in the subject. We very
briefly summarize this history:

\medskip
\noindent $\bullet$ The standard metric on $S^n$ (Gauss-Riemann).

\noindent $\bullet$ The squashed metrics on $S^{4n+3}$ (Jensen, 1973 \cite{Jen73}).

\noindent $\bullet$ A homogeneous Einstein metric on $S^{15}$ (Bourguignon and
Karcher, 1978 \cite{BuKa78}).

\noindent $\bullet$ These are all homogeneous Einstein metrics on $S^n$ and they are the
only such metrics up to homothety (Ziller, 1982 \cite{Zil82}).

\noindent $\bullet$ Infinite sequences of inhomogeneous Einstein metrics on
$S^5,S^6,S^7,S^8$ and $S^9$ (B\"ohm, 1998 \cite{Boe98}).

In contrast, we have been able to prove the following striking
results:

\begin{theorem}\cite{BGK03}
There exists inequivalent families of Einstein metrics on all odd
dimensional spheres, the number of deformation classes of \Se
structures grows double exponentially with dimension. Some of
these have moduli, the largest of which the number of moduli grows
double exponentially with dimension. There exist Einstein metrics
on all homotopy spheres in dimension $7$ and all homotopy spheres
in dimension $4n+1$ that bound parallelizable manifolds. We obtain
at least 68 inequivalent deformation classes of \Se metrics on
$S^5$, and 8610 deformation classes of \Se metrics on homotopy
spheres that are homeomorphic to $S^7.$
\end{theorem}

All the Einstein metrics in this theorem are actually \Se.

\begin{theorem}\cite{BGKT03}
There exists Einstein metrics on all homotopy spheres
in dimension $11$ and $15$ that bound parallelizable manifolds,
that is, all 992 homotopy spheres in $bP_{12}$ and all 8128
homotopy spheres in $bP_{16}$ admit \Se metrics. The distribution
of the 8610 deformation classes of \Se metrics on homotopy spheres
$\grS_i$ in dimension 7 is given by
$(n_1,\ldots,n_{28})=(376,336,260,294,231, 284,322,402,317,309,
252,\break 304,258,
390,409,352,226,260,243,309,292,452,307,298,230,307,264,353),$
where $n_i$ is the number of deformation classes of \Se metrics on
$\grS_i,$ and $\grS_1$ is the Milnor generator. For each diffeomorphism type in dimension 7 
there exist \Se metrics with moduli. The standard
$S^7=\grS_{28}$ admits an 82-dimensional family of \Se metrics.
\end{theorem}

As mentioned above B\"ohm \cite{Boe98} found infinite sequences of Einstein metrics on 
certain spheres, but it appears that our work \cite{BGK03,BGKT03} gives the first examples of 
Einstein metrics on spheres with a positive lower bound on the dimension of the moduli space. 
See \cite{Wan99} for a very nice review for obtaining Einstein metrics from symmetry 
techniques or bundle constructions. 
\bigskip

\section{Homotopy Spheres}

In 1956 John Milnor \cite{Mil56} stunned the mathematical world by
constructing smooth manifolds $\grS^7$ that are homeomorphic but
not diffeomorphic to $S^7.$ This began the field of Differential
Topology. For homotopy spheres the situation was described in a
foundational paper by Kervaire and Milnor in 1963 \cite{KeMil63}.
We briefly summarize their results. Kervaire and Milnor defined an
Abelian group $\Theta_n$ which consists of equivalence classes of
homotopy spheres of dimension $n$ that are equivalent under
oriented h-cobordism. By Smale's h-cobordism theorem this implies
equivalence under oriented diffeomorphism. The group operation on
$\Theta_n$ is connected sum. Now $\Theta_n$ has an important
subgroup $bP_{n+1}$ which consists of equivalence classes of those
homotopy spheres which are the boundary of a parallelizable
manifold. It is the subgroup $bP_{2n}$ that interests us. Kervaire
and Milnor proved:
\medskip
\begin{itemize}

\item $bP_{2m+1}=0.$\medskip

\item $bP_{4m}$ ($m\geq2$) is cyclic of order
$|bP_{4m}|=2^{2m-2}(2^{2m-1}-1)~\hbox{numerator}~\!\!\bigl({4B_m\over
m}\bigr),$ where $B_m$ is the $m$-th Bernoulli number. Thus, for
example
 $|bP_8|=28, |bP_{12}|=992, |bP_{16}|=8128, |bP_{20}|=130,816.$
\medskip
\item $bP_{4m+2}$ is either $0$ or $\bbz_2.$
\end{itemize}

\medskip
Determining which $bP_{4m+2}$  is $\{0\}$ and which is $\bbz_2$
has proven to be difficult in general, and is still not completely
understood. If $m\neq 2^i-1$ for any $i\geq 3,$ then Browder
\cite{Bro69} proved that $bP_{4m+2}=\bbz_2.$ However, 
$bP_{4m+2}$ is the identity for $m=1,3,7,15$, due to several
people \cite{MaTa67,BaJoMa84}. See \cite{La00} for a recent survey of results in this
area and complete references. The answer is still unknown in the
remaining cases. Using surgery  Kervaire was the first to show
that there is an exotic sphere in dimension $9.$ His construction
works in all dimensions of the form $4m+1$, but as just discussed
they are not always exotic.

\section{Sasakian Geometry}

What are now called {\it Sasakian structures} were first
introduced by Sasaki in 1960 \cite{Sas60} and were subsequently
developed mainly in the Japanese literature (See \cite{YK} and
\cite{Bl76a,Bl02} for complete references). They turn out to be an
odd dimensional version of K\"ahler structures.

\begin{definition-proposition}
Then a Riemannian manifold $(M,g)$ is called a {\bf
Sasakian manifold} if any one, hence all, of the following equivalent
conditions hold:
\begin{enumerate}
\item There exists a Killing vector field $\xi$ of unit length on $M$
so that the tensor field $\Phi$ of type $(1,1)$, defined by
$\Phi(X) ~=~ -\nabla_X \xi$,
satisfies the condition
$$(\nabla_X \Phi)(Y) ~=~ g(X,Y)\xi-g(\xi,Y)X$$
for any pair of vector fields $X$ and $Y$ on $M.$
\item There exists a Killing vector field $\xi$ of unit length on $M$
so that the Riemann curvature satisfies the condition
$$R(X,\xi)Y ~=~ g(\xi,Y)X-g(X,Y)\xi,$$
for any pair of vector fields $X$ and $Y$ on $M.$
\item The metric cone on $M$
$(C(M),{\bar g})= (\bbr_+\times M, \ dr^2+r^2g)$ is K\"ahler.
\end{enumerate}
\end{definition-proposition}
\medskip

We define a 1-form $\eta$ by $\eta(X)=g(\xi,X).$ Then it follows
that $\eta$ is a contact 1-form (i.e., $\eta\wedge (d\eta)^n\neq
0$) and $\xi$ is its Reeb vector field. The tensor field $\Phi$
restricted to the contact subbundle $\cald=\ker~\eta$ is an
integrable almost complex structure defining a strictly
pseudoconvex CR structure. Thus, a Sasakian structure $\cals$ is
described by tensor fields $(\xi,\eta,\Phi,g),$ which describe a
1-dimensional Riemannian foliation $\calf_\xi,$ called the
characteristic foliation, whose transverse geometry is K\"ahler.

Given a Sasakian structure $\cals=(\xi,\eta,\Phi,g)$ on $M$ there
are many Sasakian structures that can be obtained in several
different ways. First adding a basic (with respect to $\calf_\xi$)
1-form $\grz$ to $\eta$ gives a new Sasakian structure. Second
when there are non-trivial symmetries one can deform the foliation
to obtain new Sasakian structures. Although a homothety of a
Sasakian metric is no longer Sasakian, there is a ``transverse
homothety'' \cite{YK,BG05} given by $(\xi,\eta,\Phi,g)\mapsto
(a^{-1}\xi,a\eta,\Phi,ag+(a^2-a)\eta\otimes \eta)$ whose image is
a Sasakian structure. There is also a {\it conjugate Sasakian
structure} defined by $\cals^c=(-\xi,-\eta,-\Phi,g).$

We are particularly interested in
\begin{definition}
A Sasakian manifold $(M,g)$ is {\bf \Se} if the metric $g$ is also
Einstein.
\end{definition}

Any \Se metric must have positive scalar
curvature, which follows from $\Ric(X,\xi)=2n\eta(X)$.
So any complete \Se manifold must be compact with finite
fundamental group. The following is essentially due to Tanno \cite{Tan70}

\begin{lemma}\label{tanlem}
Let $\cals=(\xi,\eta,\Phi,g)$ and $\cals'=(\xi',\eta',\Phi',g)$ be two Sasakian structures sharing
the same Riemannian metric $g,$ and suppose that $(M,g)$ is not a space form with
sectional curvature equal to $1.$ Then either
\begin{enumerate}
\item $\cals'=\cals,$
\item $\cals'=\cals^c$ the conjugate Sasakian structure, or
\item $\cals$ and $\cals'$ are part of a 3-Sasakian structure.
\end{enumerate}
\end{lemma}

We refer to \cite{BG99} for the definition of a 3-Sasakian
structure. The upshot is that the types of Sasakian structures
studied here are not compatible with a 3-Sasakian structure. Lemma
\ref{tanlem} is important for the inequivalence of the Einstein
metrics we obtain. See Lemma 5.2 below.

\subsection{Relation with Algebraic Geometry}
The following two theorems generalize results of Boothby-Wang
\cite{BoWa}, Hatakeyama \cite{Hat}, and Kobayashi \cite{Kob63} to
the case where the quotient is an orbifold. A very brief
discussion of orbifolds will then follow.

\begin{theorem} \cite{BG00a}
Let $(M,g)$ be a compact quasi-regular Sasakian manifold of
dimension $2n+1$, and let $Z$ denote the space of leaves of the
characteristic foliation $\calf_\xi.$ Then
\begin{enumerate}
\item The leaf space $Z$ is a compact complex orbifold $\calz$
with a K\"ahler metric $h$ and K\"ahler form $\gro$ which defines
a class $[\gro]$ in $H^2_{orb}(\calz,\bbz)$ in such a way that
$\pi:(M,g) \ra{1.3} (\calz,h)$ is an orbifold Riemannian
submersion. The fibers of $\pi$ are totally geodesic submanifolds
of $\cals$ diffeomorphic to $S^1.$ \item The underlying complex
space $Z$ is a polarized normal projective algebraic variety with
at worst quotient singularities. \item The orbifold $\calz$ is
Fano if and only if $\Ric_g >-2.$ In this case $Z$ is simply
connected, and as an algebraic variety is uniruled with Kodaira
dimension $\kappa(Z)=-\infty.$ \item $(M,g)$ is Sasakian-Einstein
iff $(\calz,h)$ is K\"ahler-Einstein with scalar curvature \break
$4n(n+1).$
\end{enumerate}
\end{theorem}

Here $H^*_{orb}(\calz,\bbz)$ is the orbifold cohomology due to
Haefliger \cite{Hae84} (cf. \cite{BG00a} for more detail), which
rationally, but not integrally, coincides with the ordinary
cohomology $H^*(Z,\bbq).$ The important point is that
$H^2_{orb}(\calz,\bbz)$ classifies circle V-bundles (orbibundles)
over the orbifold $\calz.$ We also have an {\it Inversion Theorem}
which allows one to construct a Sasakian manifold (orbifold) from
a polarized compact K\"ahler orbifold.
\medskip

\begin{theorem} \cite{BG00a} Let $(\calz,\gro)$ be a compact
K\"ahler orbifold with $[\gro]\in H^2_{orb}(\calz,\bbz),$ and let
$\pi: M\ra{1.3} \calz$ be the $S^1$ V-bundle over $\calz$ whose
orbifold first Chern class is $[\gro].$ Suppose further that the
local uniformizing groups of the orbifold inject into $S^1,$  the
group  of the bundle, and that $[\gro]$ is a generator in
$H^2_{orb}(\calz,\bbz).$ Then $M$ is a compact simply connected
manifold which admits a Sasakian structure $\cals$ whose basic
first Chern class $c_1(\calf_\xi)$ equals $\pi^*c^{orb}_1(\calz).$
Furthermore, there is a 1-1 correspondence between compatible
K\"ahler orbifold metrics on $\calz$ in the same K\"ahler class
and homologous (with respect to the basic cohomology) Sasakian
structures on $M.$
\end{theorem}

\medskip
Here we refer to the the basic cohomology $H^*_B(\calf)$
associated with a Riemannian foliation (cf. \cite{Ton} for
details). We now have

\subsection{A Brief Review of Orbifolds}
Since in the above description the orbifold structure of $\calz$
is crucial, we give a very brief review concentrating on the
important distinctions with the geometry of manifolds or
varieties. A given algebraic variety may have many inequivalent
orbifold structures. In categorical language the topos of sheaves
on $\calz$ may be non-standard.

Orbifolds were invented by Satake \cite{Sat56,Sat57} under the
name V-manifold, and later rediscovered and renamed by Thurston
\cite{Thu79}. Since here we are only concerned with complex
orbifolds, we give the definition in this case only.

\begin{definition}
A {\bf complex orbifold} $\calz$ is a complex space $Z$ together
with a covering of charts $\calu=\{\tU_i\},$ called {\bf local
uniformizing charts}, such that the natural projections
$\varphi_i:\tU_i\ra{1.3} U_i= \bbc^n/\grG_i$ cover $Z,$  and
$\grG_i$ is a finite subgroup of $U(n)$, called a {\bf local
uniformizing group} that satisfy certain compatibility conditions.
(See \cite{Sat57,Bai57} for details).
\end{definition}
\medskip

We often write the orbifold as $\calz=(Z,\calu),$ or just $\calz.$
The local geometry on orbifolds is much the same as on manifolds,
keeping in mind that one works equivariantly on the local
uniformizing covers. However, working on these local uniformizing
covers causes certain shifts in the presence of branch divisors as
we shall see below. For complete definitions and details
concerning orbibundles (V-bundles), orbisheaves, etc., we refer to
the literature \cite{Sat57,Bai56,Bai57,MoMr03}. A fairly
comprehensive treatment will also appear in our forthcoming book
\cite{BG05}. Suffice it to say at this stage that an {\it
orbi-object} is a family of geometric objects defined on the local
uniformizing charts that satisfy the correct compatibility
conditions. For example, the {\it canonical line bundle}
$K^{orb}_{\calz}$ of the orbifold $\calz$ is a family of line
bundles, one on each chart $\tU_i$, which is the highest exterior
power of the holomorphic cotangent bundle. Actually, instead of
the canonical orbibundle, we work with {\it canonical divisors}
which by abuse of notation we denote also by $K^{orb}_{\calz}.$ In
order to understand their relation with ordinary canonical
divisors we need to consider the orbifold singular set and its
affect on divisors.

\begin{definition}
The {\bf orbifold singular set} $\grS(\calz)$ is the subset of $Z$ where
$\grG_i\neq \hbox{id}.$
\end{definition}

We shall always assume the orbifold is effective, that is that the
local uniformizing groups act effectively. It is easy to see that
the {\it orbifold regular set} $Z/\grS(\calz)$ is a dense open
subset of $Z.$ Generally, the orbifold singular set $\grS(\calz)$
differs from the usual algebro-geometric singular set. This is due
to the well known fact that quotients of a complex space under
reflections in hyperplanes are smooth. So the usual singularity
set of an algebraic variety is generally only a subset of
$\grS(\calz).$
\medskip

\begin{definition} A {\bf branch divisor} $\grD$ of an orbifold
$\calz=(Z,\calu)$ is a $\bbq$-divisor on $Z$ of the form
$$ \grD=\sum_\gra(1-\frac{1}{m_\gra})D_\gra$$
where $D_\gra$ is a Weil divisor on $Z$ that lies in the orbifold
singular locus $\grS^{orb}(Z)$, and $m_\gra$ is the $\gcd$ of the
orders of the local uniformizing groups taken over all points of
$D_\gra$, and is called the {\bf ramification
index}\index{index!ramification} of $D_\gra.$
\end{definition}
\medskip
A crucial relation for us is given by

\begin{lemma}\label{canorbidiv}
Orbifold canonical divisors $K_\calz^{orb}$ and the usual canonical divisors $K_Z$ are related 
by
$$K_\calz^{orb}\equiv \varphi^*K_Z+\sum_\gra(1-\frac{1}{m_\gra})\varphi^*D_\gra.$$
In particular $K_\calz^{orb}\equiv \varphi^*K_Z$ if and only if there are no branch divisors.
\end{lemma}

Here $\equiv$ denotes linear equivalence, and $\varphi=\sqcup \varphi_i.$
As in the usual case the orbifold first Chern class satisfies
$c_1^{orb}(\calz)= c_1(-K_\calz^{orb});$ however, generally it is only a rational class which, of
course, can differ from the first Chern class of the algebraic variety $Z.$ In particular, Fano as 
an orbifold is different than Fano as an algebraic variety. It is easy to give examples of 
non-Fano orbifold structures that lie on an algebraic variety that is Fano. It is Fano orbifold 
structures that interest us here. 

\section{Brieskorn's work}

Ten years after Milnor's famous construction of exotic differential structures on $S^7$
Brieskorn \cite{Bri66} showed how one could construct all homotopy spheres $\grS^{2n-1}\in
bP_{2n}$ explicitly.

\medskip
\noindent $\bullet$ Consider an $n+1$-tuple of positive integers $\bfa=(a_0,\dots,a_n)$
and Brieskorn-Pham polynomials in $\bbc^{n+1}$ with $a_i>1:$
$$f=z_0^{a_0}+\cdots +z_n^{a_n}.$$
\medskip

\noindent $\bullet$ Construct the link: $L(\bfa)=f^{-1}(0)\cap S^{2n+1}.$ By
the Milnor Fibration Theorem $L(\bfa)$ is $(n-2)$-connected.
\medskip

\subsection{Brieskorn Graph Theorem} To $\bfa$
one associates a graph $G(\bfa)$ whose $n+1$ vertices are labeled by $a_0,\dots,a_n.$
Two
vertices
$a_i$ and $a_j$ are connected if and only if $\gcd(a_i,a_j)>1.$ Let
$G(\bfa)_{ev}$
denote the connected component of $G(\bfa)$ determined by the even integers.
\medskip

\begin{theorem} \cite{Bri66}The link $L(\bfa)$ is homeomorphic to the
$(2n-1)$- sphere
if and only if either of the following hold:
\begin{enumerate}
\item{(i)} $G(\bfa)$ contains at least two isolated points, or
\item{(ii)} $G(\bfa)$ contains a unique
 odd isolated point and $G(\bfa)_{ev}$ has an odd number
of vertices with $\gcd(a_i,a_j)=2$
for any distinct $a_i,a_j\in G(\bfa)_{ev}$.
\end{enumerate}
\end{theorem}

\noindent One considers the cases $n=2m$ and $n=2m+1$ separately.
In particular, consider the polynomials:
$$\leqno{(1)}\qquad f=z_0^{6k-1}+z_1^3+z_2^2+\cdots +z_{2m}^2$$
and
$$\leqno{(2)}\qquad f=z_0^p+z_1^2+z_2^2+\cdots +z_{2m+1}^2 \quad (p~{\rm odd}).$$
These are examples of types (i) and (ii), respectively, of the
Graph Theorem. Brieskorn shows that with these two polynomials one
can describe all exotic spheres that bound parallelizable
manifolds! Using these polynomials we proved \cite{BGN03b} that
all such exotic spheres admit Sasakian metrics of positive Ricci
curvature, but the existence of Einstein metrics is another
matter.
\medskip

\subsection{The Diffeomorphism types}
Modulo the problems with $bP_{4m+2}$ mentioned earlier, Brieskorn
determined the diffeomorphism type of any Brieskorn manifold
satisfying the conditions of his graph theorem. Assume the
conditions of the graph theorem are satisfied, so the link is a
homotopy sphere $\grS^{2n+1}.$

\noindent $\bullet$ For $bP_{4m}$ the diffeomorphism type is
determined by the Hirzebruch signature $\grt(V_{4m})$ of the
Milnor fiber $V_{4m}$ whose boundary is $\grS^{4m-1}.$ Then
\begin{eqnarray}\label{Brieskornformula} 
\grt(V_{4m}(\bfa)) &= \#\bigl\{\bfx\in \bbz^{2m+1}
~|~0<x_i<a_i~\hbox{and}~0<\sum_{j=0}^{2m}{x_i\over a_i}
<1~{\rm mod}~ 2 \bigr\}\cr
& - \#\bigl\{\bfx\in \bbz^{2m+1}
~|~0<x_i<a_i~\hbox{and}~1<\sum_{j=0}^{2m}{x_i\over a_i} <2~{\rm mod}~ 2 \bigr\}
\end{eqnarray}

The homotopy sphere $\grS^{4m-1}_i$ is then determined by
$$i={1\over 8}\tau(V_{4m}(\Sigma_i)) {\rm mod}~ |bP_{4m}|.$$

\noindent $\bullet$ For $bP_{4m+2}\not=0$ the diffeomorphism type
is determined by the so-called {\it Arf invariant}:
$$C(V_{4m+2}(\bfa))\in \{0,1\}.$$
Then $\grS^{4m+1}$ is the standard sphere if
$C(V_{4m+2}(\bfa))=0,$ and the Kervaire sphere if
$C(V_{4m+2}(\bfa))=1.$ Furthermore, we obtain the Kervaire sphere
if and only if condition (ii) of the Brieskorn Graph Theorem holds
and the one isolated point, say $a_0,$ satisfies $a_0\equiv \pm
3~{\rm mod}~ 8.$

\section{Sasakian Geometry on Links}

In 1978 Takahashi \cite{Tak78} showed that the Brieskorn manifolds
$L(\bfa)$ naturally admit Sasakian structures. Recall that a
Riemannian manifold $(M,g)$ is Sasakian if its metric cone
$(M\times \bbr^+,dr^2+r^2g)$ is K\"ahler. \medskip

\subsection{The weighted Sasakian structure on links in $S^{2n+1}$} Let
$\bfw=(w_0,\cdots,w_n)$ be a vector whose components are positive integers.
The ``weighted'' Sasakian structure
$(\xi_\bfw,\eta_\bfw,\Phi_\bfw,g_\bfw)$ which in the standard coordinates
$\{z_j=x_j+iy_j\}_{j=0}^n$ on $\bbc^{n+1}=\bbr^{2n+2}$ is determined by
$$\eta_\bfw = {\sum_{i=0}^n(x_idy_i-y_idx_i)\over\sum_{i=0}^n
w_i(x_i^2+y_i^2)}, \qquad \xi_\bfw
=\sum_{i=0}^nw_i(x_i\partial_{y_i}-y_i\partial_{x_i}),$$ and the
standard Sasakian structure $(\xi,\eta,\Phi,g)$ on $S^{2n+1}.$
This gives the weighted Sasakian structure on $S^{2n+1}$ described
previously \cite{BG01b}.
\medskip

Now consider a Brieskorn-Pham polynomial
$$f=z_0^{a_0}+\cdots +z_n^{a_n}$$
such that $a_iw_i=d$ for each $i$.  Then the vector field
$\xi_\bfw$ restricts to the link $L(\bfa)$ and the embedding
$L(\bfa)\hookrightarrow S^{2n+1}$ induces a  weighted Sasakian
structure on $L(\bfa)$. The integer $d$ is called the {\bf
weighted degree} of $f.$

The flow of $\xi_\bfw$ gives a locally free $S^1$ action on both
$L(\bfa)$ and $S^{2n+1}.$ The quotient by this action on
$S^{2n+1}$ is a weighted projective space $\bbp(\bfw)$, and the
quotient on $L(\bfa)$ is a K\"ahler orbifold $\calz(\bfa)$ and we
have a commutative diagram
$$\begin{matrix} L(\bfa) &\ra{2.5}& S^{2n+1}_\bfw& \\
  \decdnar{\pi}&&\decdnar{} & \\
   \calz(\bfa) &\ra{2.5} &\bbp(\bfw),
\end{matrix}$$
where the horizontal arrows are Sasakian and K\"ahlerian embeddings,
respectively, and the vertical arrows are principal $S^1$ V-bundles and
orbifold Riemannian submersions.
\medskip

From the discussion above and Theorem 3.3 we have
\begin{lemma} Let $L(\bfa)$ be as above. Then the following hold:
\begin{enumerate}
\item The orbifold $\calz(\bfa)$ is Fano (i.e.,
$c^{orb}_1(\calz(\bfa))>0$) if and only if $|\bfw|:=\sum_i w_i>d.$
\item $L(\bfa)$ admits a compatible \Se metric if and only if
$\calz(\bfa)$ admits a compatible K\"ahler-Einstein orbifold
metric of scalar curvature $4n(n-1).$
\end{enumerate}
\end{lemma}

Our aim is to prove the existence of \Se metrics on the links of
Brieskorn-Pham singularities which represent homotopy spheres. By
Lemma 5.1 this is equivalent to proving the existence of \Ke
orbifold metrics on the quotient orbifolds. We also want to know
that inequivalent \Se structures imply inequivalent Einstein
metrics.
\begin{lemma}
Let $\cals=(\xi,\eta,\Phi,g)$ and $\cals'=(\xi',\eta',\Phi',g)$ be
non-conjugate Sasakian structures on $L(\bfa)$ sharing the
same Riemannian metric $g$ which is not of constant curvature.
Suppose also that $L(\bfa)$ is a homotopy sphere and satisfies
$|\bfw|-d< \frac{n}{2}{\rm min}\{w_i\}.$ Then $\cals=\cals'.$
\end{lemma}

This follows essentially from Lemma 3.2. The inequality implies
that the orbifold $\calz(\bfa)$ does not admit a holomorphic
contact structure which implies that $\cals$ and $\cals'$ cannot
be part of a 3-Sasakian structure \cite{BG97}. All of the Einstein metrics obtained in 
\cite{BGK03,BGKT03} satisfy the inequality in Lemma 5.2.
\medskip

\subsection{Existence: Solving the Monge-Amp\`ere equation}
The continuity method developed over the years by a number of
mathematicians (Aubin,Yau, Tian-Yau, Tian, Siu, Nadel,
Demailly-Koll\'ar) proves the existence of a \Ke metric on a
compact orbifold by proving an openness-closedness condition.
Openness follows from an inverse function theorem argument and, in
general, there are obstructions to the closedness. Our aim is to
find a family of functions $\phi_t$ and numbers $C_t$ for  $t\in
[0,1]$, normalized by the condition
$\int_\calz\phi_t\omega_0^{n-1}=0$, such that they satisfy the
Monge-Amp\`ere equation
$$
\log{(\omega_0+{i\over 2\pi}\partial\bar\partial\phi_t)^{n-1}\over
\omega_0^{n-1}}+t(\phi_t+f)+C_t=0.
$$
We start with $\phi_0=0, C_0=0$ and if we can reach $t=1$,
we get a \Ke metric
$$
\omega_1=\omega_0+{i\over 2\pi}\partial\bar\partial\phi_1.
$$
To find such a \Ke metric we require that the values of $t$ for
which the Monge--Amp\`ere equation is solvable approach a critical
value $t_0\in [0,1]$, a subsequence of the $\phi_t$ converges to a
function $\phi_{t_0}$ which is the sum of a $C^{\infty}$ and of a
plurisubharmonic function. As discussed by Tian \cite{Tia87} we
only need to prove that
$$
\int_\calz e^{-\gamma \phi_{t_0}}\omega_0^{n-1}<+\infty \quad
{\rm for~some}~\gamma>{n-1\over n}.
$$
Alternatively this condition can be phrased in terms of multiplier
ideal sheaves and a condition known in algebraic geometry as {\it
Kawamata log terminal (klt)}. Using Demailly and Koll\'ar
\cite{DeKo01} we get a \Ke metric on $\calz(\bfa)$, hence an \Se
metric on $L(\bfa),$ if there is a $\grg>\frac{n-1}{n}$ such that
for every weighted homogeneous polynomial $g\neq 0$ of degree
$s(|\bfw|-d)$, the function $|g|^{-\frac{\grg}{s}}$ is locally
$L^2.$ We remark that the uniform boundedness also needs to be
proven along singular orbifold divisors. In algebraic geometry
this is accomplished by the {\it inversion of adjunction}.
\medskip

Putting this together we arrive at our main operational result.

\begin{theorem}\cite{BGK03} The orbifold $\calz(\bfa)$ is Fano and
has a \Ke metric if it satisfies the condition
$$
1<\sum_{i=0}^n{1\over a_i}<
1+{n\over n-1}\min_{i,j}\Bigl\{{1\over a_i}, {1\over b_ib_j}\Bigr\}.
$$
where $b_j=\gcd(a_j,C^j)$ with $C^j={\rm lcm}(a_i:i\neq j).$
\end{theorem}
\medskip

The left hand inequality is the Fano condition while the right
hand inequality is the klt condition.

We can deform Brieskorn-Pham polynomials by adding arbitrary
monomials of weighted degree $d.$ This gives large moduli spaces
of \Ke metrics on  $\calz(\bfa)$ which is represented projectively
by the deformed weighted homogeneous polynomials.  Then one
obtains moduli spaces of \Se metrics by using Lemma 3.2. The
parameter count comes from writing down all possible monomials
satisfying the requisite conditions, and using the automorphisms
of the weighted projective spaces $\bbp(\bfw)$ to obtain normal
forms. A large lower bound on the dimension of the \Se moduli space then arises from a large 
gcd between two or more of the integers in the sequence $\bfa =(a_0,\cdots,a_n).$

\subsection{Satisfying the Inequalities} We are interested in finding sequences of integers
$\bfa=(a_0,\cdots,a_n)$ with $n\geq 3$ that satisfy the
inequalities of Theorem 5.2 as well as the conditions of the
Brieskorn Graph Theorem. For low values of $n$ the inequalities of
Theorem 5.3 are easily solved on a computer. However, it is not
only instructive to obtain some solutions by hand, but we can
prove some interesting results by judicious choices of sequences.
We consider a very important sequence that {\it does not} satisfy
the inequalities, but marks the borderline of the left hand
inequality. It is the so-called {\it extremal sequence} or {\it
Sylvester's sequence} \cite{GKP89} determined by the recursion
relation
\begin{equation*}\label{extremseq}
c_{k+1}=1+c_1\cdots c_k=
c_k^2-c_k+1
\end{equation*}
beginning with $c_1=2.$
It starts as
$$
2,3,7,43,1807, 3263443, 10650056950807,...
$$
The importance of this sequence is that it satisfies
\begin{equation*}\label{recipext}
\sum_{i=0}^n\frac1{c_i}=1-\frac{1}{c_0\cdots c_n}.
\end{equation*}
So we see that sequences of the form
$\bfa=(a_0=c_0,\dots,a_{n-1}=c_{n-1},a_n)$ satisfy the left hand
inequality as long as $a_n<c_0\cdots c_{n-1}.$ Furthermore, by
construction the first $n$ elements of such sequences are all
relatively prime to each other, so the conditions of the Brieskorn
Graph Theorem are automatically satisfied. It remains to analyze
the right hand inequality.

The troublesome part of the inequalities in Theorem 5.3 is the
computation of the $b_i$. However, for general sequences if the
$a_i$ are all pairwise relatively prime, $b_i=1$ for all $i,$ so again with
the order $a_0<a_1<\cdots <a_n$ we see that in this case the
inequalities become
$$1<\sum_{i=0}^n\frac1{a_i}< 1+\frac{n}{n-1}\frac{1}{a_n}.$$
Applying this to our special sequence
$\bfa=(a_0=c_0,\dots,a_{n-1}=c_{n-1},a_n),$ we see that the right
hand inequality automatically holds as long as $a_n$ is relatively
prime to the $c_i$'s, a condition that is easy to satisfy for all
$n.$ This gives a huge number of examples. However, since the $a_i$'s are pairwise relatively 
prime these sequences give no moduli.

Another approach which gives a large number of sequences, but now some will have moduli, 
is to notice that $b_i\leq a_i$, so it is sufficient to satisfy the
following stronger restriction:
$$
1<\sum_{i=0}^n\frac1{a_i}<
1+\frac{n}{n-1}\min_{i,j}\Bigl\{\frac1{a_ia_j}\Bigr\}=
1+\frac{n}{n-1}\cdot \frac1{a_{n-1}a_n}.
$$
By direct computation this is satisfied if $c_n-c_{n-1}<a_n<c_n$.
At least a third of these numbers are relatively prime to
$a_1=2$ and to $a_2=3$, thus we conclude

\begin{proposition}\cite{BGK03} Our methods yield at least
$\frac13(c_{n}-1)\geq \frac13(1.264)^{2^{n}}-0.5$
inequivalent families of
\Se metrics on (standard and exotic) $(2n-1)$-spheres.
\end{proposition}

If $2n-1\equiv 1\mod 4$ then all these metrics are on the standard
sphere. If $2n-1\equiv 3\mod 4$ then all these metrics are
on on both standard and exotic spheres but we cannot say anything
in general about their distribution.

\begin{example} Consider sequences of the form
$\bfa=(2,3,7,m)$. By explicit calculation, the corresponding link
$L(\bfa)$ gives a \Se metric on $S^5$ if $5\leq m\leq 41$ and $m\neq 7$
is relatively prime to at least two of $2,3,7$. This is satisfied
in $27$ cases. For example, the sequence $\bfa=(2,3,7,35)$ is
especially noteworthy. If $C(u,v)$ is any sufficiently general
homogeneous septic polynomial, then the link of
$$
x_1^2+x_2^3+C(x_3,x_4^5)
$$
also gives a \Se metric
on $S^5$.  The relevant automorphism group of $\bbc^4$ is
$$
(x_1,x_2,x_3,x_4)\mapsto
(x_1,x_2,\alpha_3x_3+\beta x_4^5,\alpha_4x_4).
$$
Hence we get a $2(8-3)=10$ real dimensional family of \Se metrics
on $S^5$. 
\end{example}

\begin{example} A similar analysis to the previous example shows, for example, that 
the sequence $\bfa=(2,3,7,43,43\cdot31)$
gives a standard 7-sphere with a $2(44-3)=82$-dimensional family
of \Se metrics on $S^7$. Similarly, the sequence $\bfa=(2,3,7,43,43\cdot39)$ gives a 
$2(44+4-5)=86$ parameter family of \Se metrics on the exotic sphere $\grS_6.$  In these 
examples we use a computer to calculate the signature of $V_8$ using Brieskorn's 
combinatorial formula \ref{Brieskornformula}.
\end{example}

The results given in Theorems 1.2 and 1.3 are obtained by computer
by inputing sequences $\bfa=(a_0,\cdots,a_n)$ and doing searches
for those satisfying the conditions of Theorems 4.1 and 5.3. We
determined the diffeomorphism type using the Brieskorn formula or
a modification of it in terms of sums of products of cotangents
due to Zagier. By Proposition 5.4 the number of deformation
classes grows double exponentially with $n,$ and it is easy to see
that the number of effective parameters as well grows double
exponentially. For example, this way one obtains 68 inequivalent
families of \Se metrics on $S^5$ (this computation can actually be
easily done without a computer). The largest family has 10
parameters. A partial computer search yielded more than $3\cdot
10^6$ cases for $S^9$ and more than $10^9$ cases for $S^{13}$,
including a 21300113901610-dimensional family. Previously, the
only known Einstein metric on $S^{13}$ was the standard one.
\medskip

Similar results for certain simply connected rational homology spheres were
obtained in \cite{BG03p}. However, the double exponential growth
is now replaced by just single exponential growth. It should should be mentioned that a 
non-vanishing second Stiefel-Whitney class $w_2$ is an obstruction for a simply connected 
manifold to admit a \Se metric \cite{BG99}. Now Smale \cite{Sm62} has classified the simply 
connected spin 5-manifolds, and the natural question arise whether all simply connected 
5-manifolds, or even all simply connected rational homology 5-spheres admit a \Se metric. 
Infinite series of rational homology 5-spheres that do admit a \Se metric were given in 
\cite{BG03p}. However, recently
Koll\'ar \cite{Kol04b} has shown that there is a torsion
obstruction to admitting even a Sasakian structure. So not all simply connected, spin, rational
homology 5-spheres can admit a \Se metric.

\section{Open Problems}

\noindent $\bullet$ {\bf Conjecture:} All homotopy spheres which
bound parallelizable manifolds admit Sasakian-Einstein metrics.
\medskip

\noindent $\bullet$ Let $K_{\rm min}$ denote the minimal value of
the sectional curvature. Can one obtain an estimate for $K_{\rm
min}$? Or better yet, a formula in terms of the weights and degree
of the BP polynomial?
\medskip

\noindent $\bullet$ Is there a bound on the dimension of the
moduli space of \Se metrics on a given manifold, (or more
generally for appropriately normalized Einstein metrics) on a
given manifold? Find formulae that depend only on dimension.
Recall that $S^{13}$ has a moduli space of \Se metrics of
dimension greater than $2.1\times 10^{13}.$
\medskip

\noindent $\bullet$ How is the moduli space of \Se metrics related
the full moduli space of Einstein metrics?
\medskip

\noindent $\bullet$ Does the moduli space of \Se metrics have an
infinite number of components? This is true for the moduli space
of deformation classes of positive Sasakian structures on spheres
\cite{BGN03b}, using work of Morita \cite{Mor75} and Ustilovsky
\cite{Ust99} on distinct contact structures. However, our proofs
of \Se metrics on spheres only yield a finite number of
deformation classes.

\medskip
\noindent $\bullet$ Does the dimension of a component depend on $K_{\rm min}?$

\def\cprime{$'$} \def\cprime{$'$} \def\cprime{$'$}
\providecommand{\bysame}{\leavevmode\hbox to3em{\hrulefill}\thinspace}
\providecommand{\MR}{\relax\ifhmode\unskip\space\fi MR }
\providecommand{\MRhref}[2]{%
  \href{http://www.ams.org/mathscinet-getitem?mr=#1}{#2}
}
\providecommand{\href}[2]{#2}


\end{document}